\documentstyle[twoside]{article}
\catcode`\@=11
\long\def\@makefntext#1{
\protect\noindent \hbox to 3.2pt {\hskip-.9pt
$^{{\eightrm\@thefnmark}}$\hfil}#1\hfill}       %CAN BE USED

\def\thefootnote{\fnsymbol{footnote}}
\def\@makefnmark{\hbox to 0pt{$^{\@thefnmark}$\hss}}    %ORIGINAL

\def\ps@myheadings{\let\@mkboth\@gobbletwo
\def\@oddhead{\hbox{}
\rightmark\hfil\eightrm\thepage}
\def\@oddfoot{}\def\@evenhead{\eightrm\thepage\hfil
\leftmark\hbox{}}\def\@evenfoot{}
\def\sectionmark##1{}\def\subsectionmark##1{}}

\oddsidemargin=\evensidemargin
\renewcommand{\thefootnote}{\fnsymbol{footnote}}

\newcounter{sectionc}\newcounter{subsectionc}\newcounter{subsubsectionc}
\renewcommand{\section}[1] {\vspace{12pt}\addtocounter{sectionc}{1}
\setcounter{subsectionc}{0}\setcounter{subsubsectionc}{0}\noindent
    {\tenbf\thesectionc. #1}\par\vspace{5pt}}
\renewcommand{\subsection}[1] {\vspace{12pt}\addtocounter{subsectionc}{1}
    \setcounter{subsubsectionc}{0}\noindent
    {\bf\thesectionc.\thesubsectionc. {\kern1pt \bfit #1}}\par\vspace{5pt}}
\renewcommand{\subsubsection}[1] {\vspace{12pt}\addtocounter{subsubsectionc}{1}
    \noindent{\tenrm\thesectionc.\thesubsectionc.\thesubsubsectionc.
    {\kern1pt \tenit #1}}\par\vspace{5pt}}
\newcommand{\nonumsection}[1] {\vspace{12pt}\noindent{\tenbf #1}
    \par\vspace{5pt}}

\newcounter{appendixc}
\newcounter{subappendixc}[appendixc]
\newcounter{subsubappendixc}[subappendixc]
\renewcommand{\thesubappendixc}{\Alph{appendixc}.\arabic{subappendixc}}
\renewcommand{\thesubsubappendixc}
    {\Alph{appendixc}.\arabic{subappendixc}.\arabic{subsubappendixc}}

\renewcommand{\appendix}[1] {\vspace{12pt}
        \refstepcounter{appendixc}
        \setcounter{figure}{0}
        \setcounter{table}{0}
        \setcounter{lemma}{0}
        \setcounter{theorem}{0}
        \setcounter{corollary}{0}
        \setcounter{definition}{0}
        \setcounter{equation}{0}
        \renewcommand{\thefigure}{\Alph{appendixc}.\arabic{figure}}
        \renewcommand{\thetable}{\Alph{appendixc}.\arabic{table}}
        \renewcommand{\theappendixc}{\Alph{appendixc}}
        \renewcommand{\thelemma}{\Alph{appendixc}.\arabic{lemma}}
        \renewcommand{\thetheorem}{\Alph{appendixc}.\arabic{theorem}}
        \renewcommand{\thedefinition}{\Alph{appendixc}.\arabic{definition}}
        \renewcommand{\thecorollary}{\Alph{appendixc}.\arabic{corollary}}
        \renewcommand{\theequation}{\Alph{appendixc}.\arabic{equation}}
        \noindent{\tenbf Appendix \theappendixc #1}\par\vspace{5pt}}
\newcommand{\subappendix}[1] {\vspace{12pt}
        \refstepcounter{subappendixc}
        \noindent{\bf Appendix \thesubappendixc. {\kern1pt \bfit #1}}
    \par\vspace{5pt}}
\newcommand{\subsubappendix}[1] {\vspace{12pt}
        \refstepcounter{subsubappendixc}
        \noindent{\rm Appendix \thesubsubappendixc. {\kern1pt \tenit #1}}
    \par\vspace{5pt}}

\newcommand{\textlineskip}{\baselineskip=13pt}
\newcommand{\smalllineskip}{\baselineskip=10pt}

\def\eightcirc{
\begin{picture}(0,0)
\put(4.4,1.8){\circle{6.5}}
\end{picture}}
\def\eightcopyright{\eightcirc\kern2.7pt\hbox{\eightrm c}}

\def\abstracts#1#2#3{{
    \centering{\begin{minipage}{4.5in}\baselineskip=10pt\footnotesize
    \parindent=0pt #1\par
    \parindent=15pt #2\par
    \parindent=15pt #3
    \end{minipage}}\par}}

\def\keywords#1{{
    \centering{\begin{minipage}{4.5in}\baselineskip=10pt\footnotesize
    {\footnotesize\it Keywords}\/: #1
    \end{minipage}}\par}}

\newcommand{\bibit}{\nineit}
\newcommand{\bibbf}{\ninebf}
\renewenvironment{thebibliography}[1]
        {\frenchspacing
     \ninerm\baselineskip=11pt
         \begin{list}{\arabic{enumi}.}
        {\usecounter{enumi}\setlength{\parsep}{0pt}
     \setlength{\leftmargin 17pt}{\rightmargin 0pt}   %FOR 10--99 ITEMS
         \setlength{\itemsep}{0pt} \settowidth
    {\labelwidth}{#1.}\sloppy}}{\end{list}}

\newcounter{itemlistc}
\newcounter{romanlistc}
\newcounter{alphlistc}
\newcounter{arabiclistc}
\newenvironment{itemlist}
        {\setcounter{itemlistc}{0}
     \begin{list}{$\bullet$}
    {\usecounter{itemlistc}
     \setlength{\parsep}{0pt}
     \setlength{\itemsep}{0pt}}}{\end{list}}

\newcommand{\fcaption}[1]{
        \refstepcounter{figure}
        \setbox\@tempboxa = \hbox{\footnotesize Fig.~\thefigure. #1}
        \ifdim \wd\@tempboxa > 5in
           {\begin{center}
        \parbox{5in}{\footnotesize\smalllineskip Fig.~\thefigure. #1}
            \end{center}}
        \else
             {\begin{center}
             {\footnotesize Fig.~\thefigure. #1}
              \end{center}}
        \fi}

\newcommand{\tcaption}[1]{
        \refstepcounter{table}
        \setbox\@tempboxa = \hbox{\footnotesize Table~\thetable. #1}
        \ifdim \wd\@tempboxa > 5in
           {\begin{center}
        \parbox{5in}{\footnotesize\smalllineskip Table~\thetable. #1}
            \end{center}}
        \else
             {\begin{center}
             {\footnotesize Table~\thetable. #1}
              \end{center}}
        \fi}

\def\@citex[#1]#2{\if@filesw\immediate\write\@auxout
    {\string\citation{#2}}\fi
\def\@citea{}\@cite{\@for\@citeb:=#2\do
    {\@citea\def\@citea{,}\@ifundefined
    {b@\@citeb}{{\bf ?}\@warning
    {Citation `\@citeb' on page \thepage \space undefined}}
    {\csname b@\@citeb\endcsname}}}{#1}}

\newif\if@cghi
\def\cite{\@cghitrue\@ifnextchar [{\@tempswatrue
    \@citex}{\@tempswafalse\@citex[]}}
\def\citelow{\@cghifalse\@ifnextchar [{\@tempswatrue
    \@citex}{\@tempswafalse\@citex[]}}
\def\@cite#1#2{{$\null^{#1}$\if@tempswa\typeout
    {IJCGA warning: optional citation argument
    ignored: `#2'} \fi}}

\def\pmb#1{\setbox0=\hbox{#1}
    \kern-.025em\copy0\kern-\wd0
    \kern.05em\copy0\kern-\wd0
    \kern-.025em\raise.0433em\box0}

\def\fnm#1{$^{\mbox{\scriptsize #1}}$}
\def\fnt#1#2{\footnotetext{\kern-.3em
    {$^{\mbox{\scriptsize #1}}$}{#2}}}

\def\fpage#1{\begingroup
\voffset=.3in
\thispagestyle{empty}\begin{table}[b]\centerline{\footnotesize #1}
    \end{table}\endgroup}

\def\runninghead#1#2{\pagestyle{myheadings}
\markboth{{\protect\footnotesize\it{\quad #1}}\hfill}
{\hfill{\protect\footnotesize\it{#2\quad}}}}
\font\tenrm=cmr10
\font\tenit=cmti10
\font\tenbf=cmbx10
\font\bfit=cmbxti10 at 10pt
\font\ninerm=cmr9
\font\nineit=cmti9
\font\ninebf=cmbx9
\font\eightrm=cmr8

\def\Z{{\mathchoice {\hbox{$\textstyle\sf Z\kern-0.4em Z$}}
{\hbox{$\textstyle\sf Z\kern-0.4em Z$}}
{\hbox{$\scriptstyle\sf Z\kern-0.3em Z$}}
{\hbox{$\scriptscriptstyle\sf Z\kern-0.2em Z$}}}}
\def\R{{\rm I\!R}}
\def\C{{\mathchoice {\setbox0=\hbox{$\displaystyle\rm C$}\hbox{\hbox
to0pt{\kern0.4\wd0\vrule height0.9\ht0\hss}\box0}}
{\setbox0=\hbox{$\textstyle\rm C$}\hbox{\hbox
to0pt{\kern0.4\wd0\vrule height0.9\ht0\hss}\box0}}
{\setbox0=\hbox{$\scriptstyle\rm C$}\hbox{\hbox
to0pt{\kern0.4\wd0\vrule height0.9\ht0\hss}\box0}}
{\setbox0=\hbox{$\scriptscriptstyle\rm C$}\hbox{\hbox
to0pt{\kern0.4\wd0\vrule height0.9\ht0\hss}\box0}}}}
\def\Q{{\mathchoice {\setbox0=\hbox{$\displaystyle\rm
Q$}\hbox{\raise
0.15\ht0\hbox to0pt{\kern0.4\wd0\vrule height0.8\ht0\hss}\box0}}
{\setbox0=\hbox{$\textstyle\rm Q$}\hbox{\raise
0.15\ht0\hbox to0pt{\kern0.4\wd0\vrule height0.8\ht0\hss}\box0}}
{\setbox0=\hbox{$\scriptstyle\rm Q$}\hbox{\raise
0.15\ht0\hbox to0pt{\kern0.4\wd0\vrule height0.7\ht0\hss}\box0}}
{\setbox0=\hbox{$\scriptscriptstyle\rm Q$}\hbox{\raise
0.15\ht0\hbox to0pt{\kern0.4\wd0\vrule height0.7\ht0\hss}\box0}}}}

\newcommand{\D}{\partial}
\newcommand{\MR}{\mathrm}
\newcommand{\W}{\wedge}
\newcommand{\ED}{'\!E}
\newcommand{\OD}{'\!O}
\newcommand{\F}[2]{\frac{#1}{#2}}
\newcommand{\LD}{\left.}
\newcommand{\RD}{\right.}
\newcommand{\LR}{\left(}
\newcommand{\RR}{\right)}
\newcommand{\LC}{\left\{}
\newcommand{\RC}{\right\}}
\newcommand{\Sum}[2]{\sum\limits_{#1=1}^{#2}}
\newcommand{\Sumo}[2]{\sum\limits_{#1}^{#2}}
\newcommand{\Tsepar}{&&&&&&&&&&&&\\[-3mm]}
\newcommand{\Overline}{\\ \hline\Tsepar}

\newcommand{\Interline}{\\ \Tsepar\hline\Tsepar}
\newcommand{\Innerline}{\\ \Tsepar}
\begin{document}

\runninghead{V. V. Kornyak}
{Computation of Cohomology of Lie Superalgebras}

\normalsize\textlineskip
\thispagestyle{empty}
\setcounter{page}{1}

%\copyrightheading{Vol. 11, No. 2 (2000) 000--000}

%\vspace*{0.88truein}

\fpage{1}
\centerline{\bf COMPUTATION OF COHOMOLOGY OF LIE SUPERALGEBRAS}
\vspace*{0.035truein}
\centerline{\bf OF VECTOR FIELDS}
\vspace*{0.37truein}
\centerline{\footnotesize V. V. KORNYAK}
\vspace*{0.015truein}
\centerline{\footnotesize\it Laboratory of Computing Technique and Automation,
Joint Institute for Nuclear Research,}
\baselineskip=10pt
\centerline{\footnotesize\it Dubna, Moscow region 141980, Russia}

\vspace*{30pt} % Kornyak addition ???
\vspace*{0.225truein}
%\publisher{(received date)}{(revised date)}

\vspace*{0.21truein}
\abstracts{The cohomology of Lie (super)algebras has many important
applications in mathematics and physics. It carries most fundamental
(``topological") information about algebra under consideration.
At present, because of the need for very tedious algebraic computation,
the explicitly computed cohomology for different classes of
Lie (super)algebras is known only in a few cases. That is why
application of computer algebra methods is important for this problem.
We describe here an algorithm and its C implementation for computing
the cohomology of Lie algebras and superalgebras.
The program can proceed finite-dimensional algebras and infinite-dimensional
graded algebras with finite-dimensional homogeneous components.
Among the last algebras  Lie algebras and superalgebras
of formal vector fields are most important.
We present some results of computation of cohomology for Lie superalgebras
of Buttin vector fields and related algebras.
These algebras being super-analogs of Poisson and Hamiltonian algebras
have found many applications to modern supersymmetric models of theoretical
and mathematical physics.
}{}{}

\vspace*{10pt}
\keywords{Cohomology, Lie superalgebras, Antibracket,
Buttin superalgebra, Symbolic computation.}

\vspace*{10pt} % Kornyak addition ???

\vspace*{1pt}\textlineskip  %) USE THIS MEASUREMENT WHEN THERE IS
\section{Introduction and basic definitions} %) SECTION HEADING
\vspace*{-0.5pt}
\noindent
There are many applications of the Lie (super)algebra
cohomology in mathematics: characteristic classes of foliations;
invariant differential operators; MacDonald-type combinatorial
identities, etc. (see\cite{Fuks} for details).
The use of cohomology in theoretical and mathematical
physics can be illustrated by the following applications:
\begin{itemlist}
 \item construction of the central ex\-ten\-si\-ons and
  deformations for Lie super\-al\-gebras;
 \item construction of the supergravity equations for $N$-extended Minkowski
  super\-spa\-ces and search for possible models for these super\-spa\-ces;
 \item study of stability for nonholonomic systems like ballbearings,
  gyroscopes, elec\-t\-ro-mechanical devices, waves in plasma, etc.;
 \item description of an analogue of the curvature tensor for nonlinear
  nonholonomic constraints;\cite{GL}
 \item new methods for the study of integrabi\-li\-ty of dynamical systems;
\pagebreak

\textheight=7.8truein       %2ND \TEXTHEIGHT IS ONLY FOR SUBSEQUENT PAGES
\setcounter{footnote}{0}
\renewcommand{\thefootnote}{\alph{footnote}}

 \item construction of so-called {\em higher order Lie algebras}\cite{Fil}
  which allow in turn to construct the {\em Nambu mechanics}\cite{Takh}
  generalizing the ordinary Hamiltonian mechanics;
 \item construction of possible invariant effective actions of
  Wess-Zumino-Witten type and the study of anomalies.\cite{D'Hoker}
\end{itemlist}

General definitions and properties of cohomology of Lie algebras
and superalgebras are described in book.\cite{Fuks}
Let us recall briefly some basic definitions.

A {\em Lie superalgebra\/} is a ${\Z}_2$-graded algebra over a commutative
ring $K$ with a unit:
$$L = L_{\bar{0}} \oplus L_{\bar{1}}, \
u \in L_\alpha, \  v \in L_\beta, \ \alpha, \beta \in {\Z}_2 =
\{\bar{0}, \bar{1}\}\Longrightarrow [u,v] \in L_{\alpha + \beta}$$
The elements of $L_{\bar{0}}$ and $L_{\bar{1}}$ are called {\em even\/}
and {\em odd}, respectively.
By definition, the {\em Lie product\/}
 (shortly, {\em bracket\/})
 $[\cdot ,\cdot ]$ satisfies the following axioms
\begin{eqnarray*}
& [u,v] = -(-1)^{p(u)p(v)} [v,u], &  {\qquad skew-symmetry,} \label{ss_s} \\
& [u,[v,w]] = [[u,v],w] +(-1)^{p(u)p(v)} [v,[u,w]], &
{\qquad \mbox{Jacobi\ identity},}
\end{eqnarray*}
where $p(a)$ is the parity of element $a \in L_{p(a)}$.
We shall assume that $K$ is a field. If $K$ is a field of characteristic
$2$ or $3$, extra axioms are needed:
$[u,u]=0$ for even $u$ in characteristic $2$ and $[v,[v,v]]=0$ for odd
$v$ in characteristic $3$.
To provide connection with enveloping algebra, characteristic $2$
requires also the existence of a {\em quadratic operator} $q$ mapping
odd elements of the algebra into even ones such that
\begin{eqnarray*}
&& q(\alpha u)=\alpha^2 q(u), \\
&& [u,v]=q(u+v)-q(u)-q(v), \\
&& [u,[u,v]]=[q(u),v],
\end{eqnarray*}
where $\alpha \in K$ and $u,v$ are odd.

A {\em module\/} over a Lie superalgebra $A$ is a vector space $M$
(over the same field $K$) with a mapping $A \times M \to M$, such that
$[a_1, a_2] m = a_1 (a_2 m) - (-1)^{p(a_1)p(a_2)} a_2 (a_1 m)$,
where $a_1, a_2 \in A$, $m \in M$.
The most important (and easy for program implementation) are {\em trivial\/}
($M$ is arbitrary vector space, e.g., $M = K$; $am = 0$), {\em adjoint\/}
($M = A; am = [a,m]$) and {\em coadjoint\/}  ($M = A'; am = \{a,m\}$
is coadjoint action) modules.

A {\em cochain complex\/} is a sequence of linear spaces $C^k$
with linear mappings $d^k$
\begin{equation}
0\to C^0\stackrel{d^0}{\longrightarrow}\cdots
\stackrel{d^{k-2}}{\longrightarrow}C^{k-1}\stackrel{d^{k-1}}
{\longrightarrow} C^k\stackrel{d^k}{\longrightarrow}C^{k+1}
\stackrel{d^{k+1}}{\longrightarrow}\cdots,
\label{cocomplex}
\end{equation}
where the linear space $C^k = C^k(A;M)$ is a super skew-symmetric
$k$-linear mapping $A \times \cdots \times A \to M$, $C^0 = M$ by definition.
The super skew-symmetry means symmetry w.r.t. transpositions of odd adjacent
elements of $A$ and antisymmetry for all other transpositions of adjacent elements.
Elements of $C^k$ are called {\em cochains}.

The linear mapping $d^k$ (or, briefly, $d$) is called a {\em differential\/}
and satisfies the following property: $d^k \circ d^{k-1} = 0$ (or $d^2 = 0$).

The cochains mapped into zero by the differential are called
{\em cocycles}, i.e., the space of cocycles is
$$Z^k = {\rm Ker}\ d^k = \{C^k \; | \; dC^k = 0\}.$$

The cochains which can be represented as differentials
of other cochains are called {\em coboundaries}, i.e., the space
of coboundaries is
$$B^k = {\rm Im}\ d^{k-1} = \{C^k \; | \; C^k = dC^{k-1}\}.$$
Any coboundary is, obviously, a cocycle.

The non-trivial cocycles, i.e., those which are not coboundaries,
form the {\em cohomology}. In other words, the cohomology is
the quotient space $$H^k(A;M) = Z^k/B^k.$$
As is seen from the initial part of cochain complex (\ref{cocomplex})
the basis elements of module $M$ can be considered as non-trivial
0-cocycles.

The explicit form of the differential for a Lie superalgebra is
\begin{eqnarray*}
& dC(e_0,\ldots,e_q;O_{q+1},\ldots,O_k)\ = & \\
&  \Sumo{i<j}{q}(-1)^j C(e_0,\ldots,e_{i-1},[e_i,e_j],\ldots,
    \widehat{e_j},\ldots,e_q;O_{q+1},\ldots,O_k)\ + &\\
& (-1)^{q+1}\Sumo{i=0}{q}\Sumo{j=q+1}{k}C(e_0,\ldots,e_{i-1},[e_i,O_j],\ldots,
   e_q;O_{q+1},\ldots,\widehat{O_j},\ldots, O_k)\ + &\\
& (-1)^{i+1}\Sumo{i=q+1}{k-1}\Sumo{j=q+2}{k}C(e_0,\ldots,e_q;O_{q+1},\ldots,
   O_{i-1},[O_i,O_j],\ldots,\widehat{O_j},\ldots, O_k) &\\
& + \Sumo{i=0}{q}(-1)^{i+1} e_i
   C(e_0,\ldots,\widehat{e_i},\ldots,e_q;O_{q+1},\ldots,O_k) &\\
& +\ (-1)^q \Sumo{i=q+1}{k} O_i
   C(e_0,\ldots,e_q;O_{q+1},\ldots,\widehat{O_i},\ldots,O_k). &
\end{eqnarray*}
Here $e_i$ and $O_i$ are even and odd elements of the algebra, respectively,
and the hat ``$\;\widehat{\ }\;$" marks the omitted elements.

Here are some properties and statements we use in the sequel.

An algebra and a module are called {\em graded\/} if they can be
presented as sums of homogeneous components in a way compatible
with the algebra bracket and the action of the algebra on the module:
$$A = \oplus_{g \in G}\, A_g, \; M = \oplus_{g \in G}\, M_g, \;
[A_{g_1}, A_{g_2}] \subset A_{g_1+g_2}, \;
A_{g_1}M_{g_2} \subset M_{g_1+g_2},$$
where $G$ is some abelian (semi)group. We assume $G = \Z$ in this paper.
To avoid confusion, we use in the sequel the terms {\em grade} and
{\em degree}
for element of $G$ and number of cochain arguments, respectively.
The grading in the algebra and module induces a grading on cochains
and, hence, in the cohomology:
$$C^{*}(A;M) = \oplus_{g \in G}\, C_g^{*}(A;M), \quad
 H^{*}(A;M) = \oplus_{g \in G}\, H_g^{*}(A;M).$$
This property allows one to compute the cohomology separately
for different ho\-mo\-ge\-ne\-ous components; this is especially useful
when the homogeneous components are finite-dimensional.

If there is an element $a_0 \in A$, such that eigenvectors
(with the same eigenvalues for a given grade) of the operator
$a \mapsto [a_0,a]$ form a (topological) basis of algebra $A$,
then $H^{*}(A) \simeq H_0^{*}(A).$ In other words, all the non-trivial
cocycles of the cohomology in the trivial module lie in the zero grade
component.
The element $a_0$ is called an {\em internal grading element}.
If eigenvectors of the operator $m \mapsto a_0m$
form also a topological basis of module $M$, then the same statement holds
for the cohomology in the module $M$: $H^{*}(A;M) \simeq H_0^{*}(A;M).$

In the case of trivial module, the exterior multiplication of cochains
provides the cohomology with a structure of graded ring, i.e.,
if $C^k$ and $C^m$ are cocycles then $C^{k+m} = C^k\W C^m$ is also
a cocycle.

If algebra $A$ contains a {\em central element\/} $Z$, i.e.,
$[Z,a]=0$ for any $a\in A$, then cochain
$C^{k+1}(a_1,\ldots,a_k, Z) = C^k(a_1,\ldots,a_k)\W C^1(Z)$
is a cocycle provided that $C^k$ is a cocycle: $dC^k=0\Rightarrow dC^{k+1}=0$
because $dC(a,Z)\sim C([Z,a]) = 0$ and the differential $d$ acts on a product
of cochains as a (super)differentiation. Due to this fact any cocycle
$C^k(a_1,\ldots,a_k)$ for algebras with {\em odd\/} center leads to an
infinite set of cocycles of the form
$C^{k+m}(a_1,\ldots,a_k,\underbrace{Z,\ldots,Z}_m\,).$ We shall
encounter such situation later in our computations.

There are also another multiplicative structures in the cohomology theory,
but we shall not use them in this work.

\section{Lie superalgebras of vector fields}
\noindent
Below a list of the main Lie superalgebras of formal vector fields is
given.\cite{Leites}
We consider some sets of even ($x_i, q_i, p_i, t$) and odd (called also
{\em Grassmann\/}) variables ($X_i, T$).
In many cases the vector fields can be expressed in terms
of generating functions. The coordinates of vector fields and generating
functions are assumed to be formal power series in the even and odd variables.
Note that all the algebras depending only on the odd variables are
finite-dimensional.
All these algebras are graded due to a prescribed grading of the variables.
There are some {\em standard\/} gradings for the variables:
all variables $x_i, q_i, p_i, X_i$ have grade 1 and the separate
variables $t, T$ have grade 2. {\em Non-standard\/} gradings
with zero or negative grades for some odd variables are possible
(and useful, as we show below) too.
The divergence-free algebras are called {\em special}.
The symbol  ${\bf Z}$ denotes the 1-dimensional center of an algebra
consisting of constants in terms of generating function, i. e.,
${\bf Z}=\MR{Span}(\{1\})$.
\begin{enumerate}
%------------------------------------------------------------
 \item {\em General vectorial superalgebra}
  ${\bf W(n \mid m)}$ or ${\bf vect(n \mid m)}$ \\
 Variables: $x_1,\ldots,x_n; X_1,\ldots,X_m$ \\
 The bracket denotes the supercommutator of vector fields of the form \\
 $\Sum{i}{n} f_i\F{\D}{\D x_i} + \Sum{k}{m} g_k\F{\D}{\D X_k}$

%------------------------------------------------------------
  \item {\em Special vectorial superalgebra}
 ${\bf S(n \mid m)}$ or ${\bf svect(n \mid m)}$
 consists of the elements
 from ${\bf W(n \mid m)}$ satisfying the divergence-free condition \\
 $\Sum{i}{n}\F{\D f_i}{\D x_i}+\Sum{k}{m}(-1)^{p(g_k)}\F{\D g_k}{\D X_k}=0$

%------------------------------------------------------------
 \item {\em Poisson superalgebra}
 ${\bf Po(2n \mid m)}$\\
 Variables: $p_1,\ldots,p_n, q_1,\ldots,q_n; X_1,\ldots,X_m$\\
 Bracket:
  $\{f,g\}_{Pb}=\Sum{i}{n}\LR\F{\D f}{\D p_i}\F{\D g}{\D q_i}
                            -\F{\D f}{\D q_i}\F{\D g}{\D p_i}\RR
        -(-1)^{p(f)}\Sum{k}{m}\F{\D f}{\D X_k}\F{\D g}{\D X_k}$\\
 Sometimes it is more convenient to redenote the odd variables
 $X_k$ and set
 $P_k = \F{1}{2}(X_k-X_{r+k}),  Q_k = \F{1}{2}(X_k+X_{r+k})$
 for $k\leq r =[m/2],$\\
 ($U = X_{2r+1}$ for $m$ odd)
 and the last sum in the bracket takes the form
 $\Sum{k}{r}(\F{\D f}{\D P_k}\F{\D g}{\D Q_k}+\F{\D f}{\D Q_k}\F{\D g}{\D
 P_k})$ (one should add the term $\F{\D f}{\D U}\F{\D g}{\D U}$ for $m$ odd).\\
 {\em Hamiltonian superalgebra\/} is
 ${\bf H(2n\mid m) = Po(2n\mid m)/Z}$\\
 {\em Special Hamiltonian superalgebra\/} ${\bf SH(0\mid m)}$ is a simple
 ideal of codimension one in ${\bf H(0\mid m)}$

%------------------------------------------------------------
 \item {\em Contact superalgebra}
 ${\bf K(2n+1 \mid m)}$\\
 Variables: $t, p_1,\ldots,p_n, q_1,\ldots,q_n; X_1,\ldots,X_m$\\
 Bracket: $\{f,g\}_{Kb} = \delta(f)\F{\D g}{\D t}-\F{\D f}{\D t}\delta(g)
                                                  -\{f,g\}_{Pb}$\\
 $\delta(f) = 2f-E(f), \
 E = \Sum{i}{n}\LR p_i\F{\D}{\D p_i} + q_i\F{\D}{\D q_i}\RR
                              + \Sum{k}{m} X_k\F{\D}{\D X_k}$

%------------------------------------------------------------
 \item {\em Buttin superalgebra\/} ${\bf B(n)}$\\
 Variables: $x_1,\ldots,x_n; X_1,\ldots,X_n$\\
 Bracket:
 \begin{equation}
 \{f,g\}_{Bb}=\Sum{i}{n}\LR\F{\D f}{\D x_i}\F{\D g}{\D X_i}
                             +(-1)^{p(f)}\F{\D f}{\D X_i}\F{\D g}{\D x_i}\RR
 \label{Bb}
 \end{equation}
 {\em Leites superalgebra\/} is ${\bf Le(n) = B(n)/Z}$

%------------------------------------------------------------
 \item {\em Special Buttin superalgebra\/} ${\bf SB(n)}$
 is subalgebra of ${\bf B(n)}$ \\ subject to the constraint
 $\Delta f = 0$ for generating function,
 where
 \begin{equation}
  \Delta = \Sum{i}{n}\F{\D^2}{\D x_i\D X_i}.
 \label{Delta}
 \end{equation}
 {\em Special Leites superalgebra\/} is ${\bf SLe(n) = SB(n)/Z}$

%------------------------------------------------------------
 \item {\em Odd contact superalgebra\/} ${\bf M(n)}$\\
 Variables: $x_1,\ldots,x_n; T, X_1,\ldots,X_n$\\
 Bracket: $\{f,g\}_{Mb}=\delta(f)\F{\D g}{\D T}
            +(-1)^{p(f)}\F{\D f}{\D T}\delta(g) - \{f,g\}_{Bb}$\\
 $\delta(f) = 2f-E(f), \
 E = \Sum{i}{n}\LR x_i\F{\D}{\D x_i} + X_i\F{\D}{\D X_i}\RR$

%------------------------------------------------------------
 \item {\em Special odd contact superalgebra\/} ${\bf SM(n)}$
  is subalgebra of ${\bf M(n)}$ subject to the constraint
  $(1-E)\F{\D f}{\D T} - \Delta f=0$
  for generating function.
\end{enumerate}
\section{Outline of algorithm and its implementation}
\noindent
To compute the cohomology one needs to solve the equation
\begin{equation}
dC^k = 0,
\label{dck}
\end{equation}
and throw away those solutions of (\ref{dck}) which can be
expressed in the form
$$C^k = dC^{k-1}.$$
In the case of finite-dimensional Lie superalgebras determining
equation (\ref{dck}) is a system of finite-dimensional homogeneous
linear algebraic equations.
In the case of infinite-dimensional graded Lie superalgebras, such as
Lie superalgebras of vector fields with even variables,
equation (\ref{dck}) is a system of linear homogeneous functional equations
with integer arguments. Unfortunately there is no general method for solving
such systems in closed form, though in a few exceptional cases such solutions
are known.

If the grading leads to finite-dimensional space of cochains
in a given grade, one can proceed just as in the case of finite-dimensional
algebra. Unfortunately the very important case of computation of cohomology
in adjoint module for infinite-dimensional algebras can not be reduced
to the set of finite-dimensional tasks at any choice of grading: in the case
of adjoint module the cochains contain both elements of algebra and dual
elements,
these elements inevitably should have opposite grades. Nevertheless,
there are important problems
(such as the Spencer cohomology playing an essential role
in the formal theory of differential equations\cite{Pommaret})
requiring computation of cohomology in adjoint module with respect
to finite-dimensional subalgebras of infinite-dimensional algebras.

There are several packages for computing
cohomology of Lie algebras and superalgebras written in {\em Reduce}
\cite{LP,PH} and {\em Mathematica}.\cite{GL} Some new results
were obtained completely or partially with the help of these
packages.\fnm{a}\fnt{a}{In particular, D. Leites informed us that
A. Shapovalov discovered one of the cocycles from  the cohomology of
special Hamiltonian superalgebra $\MR{SH}(0|4)$ with the help of
the {\em Mathematica} program\cite{GL} written by P. Grozman.
This cocycle was missed by D. Leites and D. Fuchs when they investigated
this cohomology\cite{FL} by hand.}
However, abilities of these packages are restricted by rather small problems.
We wrote a more advanced program\cite{KornProg} which allows us to consider
more difficult and real problems.

The C code of the program, of total length near 14200 lines, contains
about 400 functions realizing top level algorithms, simplification of
indexed objects, working with Grassmannian objects, exterior calculus,
linear algebra, substitutions, list processing, input and output, etc.

All operations with scalar coefficients, including input and output,
are localized in 17 functions which formats do not
depend on the nature of field of scalars $K$. Reading from the input which
field should be used, the program assigns the suitable function addresses
to the corresponding function pointer variables. This feature of the
C language allows to carry out computations over arbitrary fields
without recompiling and any loss of efficiency. Up to now we have implemented
rational numbers of arbitrary precision, i.e., the field $\Q$, its complex
extension $\Q[i]$,\fnm{b}\fnt{b}{Note that the fields $\R$ and $\C$
being non-constructive objects do not admit a computer implementation
at all.} rational functions of arbitrary parameters
(for classification problems) and the fields $\Z_p$.\fnm{c}\fnt{c}{For
efficiency reasons the prime $p$ should not exceed $46337$ on 32bit
and $3037000493$ on 64bit computers.}
Of course, other fields can be easily added if necessary.

We represent Grassmann monomials by integer numbers using one-to-one
correspondence between (binary codes of) non-negative integers and
Grassmann monomials. This representation allows one efficiently to implement
the operations with Grassmann monomials by means of the basic computer
commands.

The program performs sequentially the following steps:
\begin{enumerate}
\item {\em Reading input information.}
\item {\em Constructing a basis\/} for the algebra. The basis can be read
 from the input file; otherwise the program constructs it from
 the definition of the algebra.
 Non-trivial computations
 at this step arise only in the case of divergence-free algebras.
 The basis elements of such algebras should satisfy some conditions.
 In fact, we should construct the basis elements of the subspace given by a
 system of linear equations. The task is thereby reduced to some problem of linear
 algebra combined with shifts of indices. For example, among the divergence-free
 conditions for the special Buttin algebra ${\MR{SB}(3)}$ there are the following
 two equations
 $$ia_{ijk; XY} -(k+1)a_{i-1,j,k+1; YZ} = 0,$$
 $$ia_{ijk; XZ} +(j+1)a_{i-1,j+1,k; YZ} = 0.$$
 Here $a_{ijk; XY},\ldots$ are coefficients at the monomials $x^iy^jz^kXY,\ldots$
 in the generating function; $x, y, z$ and $X, Y, Z$ are even and odd variables,
 respectively. First of all, we have to shift indices $j$ and $k$ in the second
 equation to reduce the last terms of both equations to the same multiindices.
 Then, using some simple tricks of linear algebra, we can easily construct the
 corresponding basis element
 $$O^1_{ijk}=kx^iy^jz^{k-1}XY-jx^iy^{j-1}z^kXZ+ix^{i-1}y^jz^kYZ.$$
 As a result for ${\MR{SB}(3)}$ we have the basis:
\begin{eqnarray*}
(1)~ E^1~~&=&X Y Z\\
(2)~ E^2_{ijk}&=&kx^{i} y^{j} z^{k-1}X-ix^{i-1} y^{j} z^{k}Z\\
(3)~ E^3_{ijk}&=&kx^{i} y^{j} z^{k-1}Y-jx^{i} y^{j-1} z^{k}Z\\
(4)~ O^1_{ijk}&=&kx^{i} y^{j} z^{k-1}X Y-jx^{i} y^{j-1} z^{k}X Z+ix^{i-1} y^{j} z^{k}Y Z\\
(5)~ O^2_{ijk}&=&x^{i} y^{j} z^{k}
\end{eqnarray*}
\item {\em Constructing the commutator table\/} for the algebra
 (if this table has not been read from the input file).\\
We illustrate this step by non-zero commutators of $\MR{SB}(3)$ generated by the program:
$$\begin{array}{llll}
(1)&[E^2_{ijk}, E^2_{lmn}]&=&(ni-lk)E^2_{i+l-1,j+m,k+n-1}\\
(2)&[E^2_{ijk}, E^3_{lmn}]&=&\F{nkj-mk^2+mk}{n+k-1}E^2_{i+l,j+m-1,k+n-1}\\
& & &+\F{n^2i-nlk-ni}{n+k-1}E^3_{i+l-1,j+m,k+n-1}\\
(3)&[E^3_{ijk}, E^3_{lmn}]&=&(nj-mk)E^3_{i+l,j+m-1,k+n-1}\\
(4)&[E^2_{ijk}, O^1_{lmn}]&=&(ni-lk)O^1_{i+l-1,j+m,k+n-1}\\
(5)&[E^3_{ijk}, O^1_{lmn}]&=&(nj-mk)O^1_{i+l,j+m-1,k+n-1}\\
(6)&[E^1, O^2_{ijk}]&=&-O^1_{ijk}\\
(7)&[E^2_{ijk}, O^2_{lmn}]&=&(ni-lk)O^2_{i+l-1,j+m,k+n-1}\\
(8)&[E^3_{ijk}, O^2_{lmn}]&=&(nj-mk)O^2_{i+l,j+m-1,k+n-1}\\
(9)&[O^1_{ijk}, O^2_{lmn}]&=&\F{nj-mk}{n+k}E^2_{i+l,j+m-1,k+n}
+\F{-ni+lk}{n+k}E^3_{i+l-1,j+m,k+n}\\
\end{array}$$
\item\label{stepgenform}
 {\em Creating the general form\/} of expressions for coboundaries
 and determining  equations for cocycles.
\item {\em Transition to a particular grade\/} in general expressions.
 At this step expressions for coboundaries take the form
 ${\bf x = bt}$, equations for cocycles take the form ${\bf Zx = 0}$,
 where vector ${\bf x}$ corresponds to $C^k$, parameter vector ${\bf t}$
 corresponds to $C^{k-1}$, matrices ${\bf Z, b}$ correspond to the differential
 $d$. All these vector spaces are finite-dimensional for any particular grade.
\item {\em Computing the quotient space $H^k(A;M) = Z^k/B^k$}.
 Here the cocycle subspace $Z^k$ is given by relations ${\bf Zx = 0}$, and
 the coboundary subspace $B^k$ is given parametrically by ${\bf x = bt}$.

 Substeps:
 \begin{enumerate}
 \item Eliminate ${\bf t}$ from  ${\bf x = bt}$ to get equations
  ${\bf Bx = 0}$
 \item Reduce both relations ${\bf Bx = 0}$ and ${\bf Zx = 0}$ to the
  canonical (row echelon) form by Gauss elimination.
  If ${\rm rank} {\bf B} = {\rm rank} {\bf Z}$,
  then there is no non-trivial cocycle;
  otherwise go to Substep (\ref{substepc}).
 \item\label{substepc}
  Set ${\bf Bx = y}$ and substitute these relations into ${\bf Zx = 0}$
  to get relations ${\bf Ay = 0}$. The {\em parametric} (non-leading)
  $y'$s of the last relations are non-trivial cocycles; that is,
  they form a basis of the cohomology.
 \end{enumerate}
 In fact, the above procedure is based on the relation for quotient spaces
 $$Z/B = \F{Y/B}{Y/Z},$$
 where $Y$ is an artificially introduced space,
 combining the above $x'$s and $y'$s.
 \item\label{stepoutput}
 {\em Output the non-trivial cocycles.}
  The program can output results in 2D ASCII, \LaTeX\
  and standard for usual computer algebra systems 1D forms.
  The last form of output is useful for investigating the structure
  of cohomology ring with the help of the systems like
  {\em Maple}, {\em Mathematica} or {\em Reduce}. The operations
  in such investigations (multiplications and comparisons of cochains)
  are not difficult from computational point of view, but
  interactive abilities of the above systems are very convenient
  for the analysis of the cohomology ring.
\end{enumerate}
To split the whole task to smaller ones
Steps (\ref{stepgenform}--\ref{stepoutput}) are executed separately
for even and odd parts of the cochain complex.

\subsection{Example of output file: Computation of $H^5_0(\MR{SLe}(2))$}
\noindent
The below output demonstrates computation of 5-cocycle in grade 0 from
cohomology in trivial module for special Leites superalgebra of vector fields
in superdimension $(2|2)$. Here are some explanations to this output.
The brackets $<\cdots>$ include a comment in input file. Some elements of input
are optional. If they are omitted, the program either constructs some standard
ones or asks to input them from keyboard. $g(a)$ is a $\Z$-grade of element $a$.
$E_{ij}$ and $O^1,\ O^2_{ij}$ are even and odd basis elements of superalgebra.
$'\!E_{ij},\ '\!O^1,\ '\!O^2_{ij}$ are dual elements to corresponding basis elements.
The vertical and horizontal dots mean that we have omitted for brevity some long
(sub)expressions in this illustrative example.
Output of determining equations for cocycles may be useful. Sometimes
(for a low degree cohomology) one can see the general solution
for these equations.
Almost all elements of output (excepting the resulting cocycle)
can be suppressed by corresponding settings in the initiating file.
$t_i$ are arbitrary parameters describing the space of coboundaries.
This computation gives one basis element of the cohomology, but the program
produces also its four equivalent forms: one can choose any of them (or their
linear combination)
in order to get more compact or symmetrical expression.
Note that in the case of several basis elements of cohomology,
these alternative forms may be linear combinations of the original ones.
\newpage
\begin{verbatim}
Input file: D:\Kornyak\LieCohomology\In\test.in
Input data:
<* Special Leites superalgebra SLe(n) = SB(n)/Z *>
Even variables: x; y. < Optional >
Grading for even variables: 1; 1. < Optional >
Odd variables:  X; Y. < Optional >
Grading for odd variables: -1; -1. < Optional >
Module type: Trivial. < Coadjoint Adjoint>
Special Leites superalgebra: 2.
Cohomology number: 5. < Optional >
Grade: 0. < Optional >
\end{verbatim}
\noindent
Even variables of vector field:$\ x\ y$;$\ \; g(x)=1\; g(y)=1.$\\
\noindent
Odd variables of vector field: $\ X\ Y$;$\ \; g(X)=-1\; g(Y)=-1.$\\
\noindent
Basis elements of Lie superalgebra:\label{outputexample}
$$\begin{array}{llllll}
(1)&E_{ij}&=&jx^{i} y^{j-1}X-ix^{i-1} y^{j}Y;&g(E_{ij})=i+j-2;&i \geq 0,\ j \geq 0.\\
(2)&O^1&=&X Y;&g(O^1)=-2&\\
(3)&O^2_{ij}&=&x^{i} y^{j};&g(O^2_{ij})=i+j;&i \geq 0,\ j \geq 0.
\end{array}$$

\noindent
Non-zero commutators of Lie superalgebra:
$$\begin{array}{llll}
(1)&[E_{ij}, E_{kl}]&=&(li-kj)E_{i+k-1,j+l-1}\\
(2)&[E_{ij}, O^2_{kl}]&=&(li-kj)O^2_{i+k-1,j+l-1}\\
(3)&[O^1, O^2_{ij}]&=&-E_{ij}\\
\end{array}$$

\noindent
Expression for even coboundaries:
\begin{eqnarray*}
dC_{\bar{0}}^4&=&\LC(ro-qp)C\LR E_{ij},E_{kl},E_{mn},E_{o+q-1,p+r-1}\RR\RD+\cdots\\
&&\LD\cdots+(-ri+qj)C\LR O^2_{i+q-1,j+r-1},O^2_{kl},O^2_{mn},O^2_{op}\RR\RC\\
&&\ED_{ij}\W\OD^2_{kl}\W\OD^2_{mn}\W\OD^2_{op}\W\OD^2_{qr}
\end{eqnarray*}

\noindent
Expression for odd coboundaries:
\begin{eqnarray*}
dC_{\bar{1}}^4&=&\LC(-pm+on)C\LR E_{ij},E_{kl},E_{m+o-1,n+p-1},O^1\RR\RD+\cdots\\
&&\LD\cdots+C\LR E_{op},O^2_{ij},O^2_{kl},O^2_{mn}\RR\RC
\OD^1\W\OD^2_{ij}\W\OD^2_{kl}\W\OD^2_{mn}\W\OD^2_{op}
\end{eqnarray*}

\noindent
Determining equation for even cocycles:
\begin{eqnarray*}
dC_{\bar{0}}^5&=&\LC(-tq+sr)C\LR E_{ij},E_{kl},E_{mn},E_{op},E_{q+s-1,r+t-1}\RR\RD+\cdots\\
&&\LD\cdots+C\LR E_{qr},O^2_{ij},O^2_{kl},O^2_{mn},O^2_{op}\RR\RC
\OD^1\W\OD^2_{ij}\W\OD^2_{kl}\W\OD^2_{mn}\W\OD^2_{op}\W\OD^2_{qr}
\end{eqnarray*}

\noindent
Determining equation for odd cocycles:
\begin{eqnarray*}
dC_{\bar{1}}^5&=&\LC(ro-qp)C\LR E_{ij},E_{kl},E_{mn},
E_{o+q-1,p+r-1},O^1\RR\RD+\cdots\\
&&\LD\cdots+(-ti+sj)C\LR O^2_{i+s-1,j+t-1},O^2_{kl},
O^2_{mn},O^2_{op},O^2_{qr}\RR\RC\\
&&\ED_{ij}\W\OD^2_{kl}\W\OD^2_{mn}\W
\OD^2_{op}\W\OD^2_{qr}\W\OD^2_{st}
\end{eqnarray*}

\noindent
Coboundary component expressions in grade 0:
$$\begin{array}{ccc}
C\LR E_{01},E_{02},E_{03},E_{10},E_{12}\RR&=&2t_1-3t_3\\
\vdots&\vdots&\vdots\\
C\LR E_{10},O^1,O^2_{10},O^2_{10},O^2_{10}\RR&=&0
\end{array}$$
where
$$\begin{array}{ccc}
t_1&=&C\LR E_{01},E_{02},E_{03},E_{11}\RR\\
\vdots&\vdots&\vdots\\
t_{246}&=&C\LR O^1,O^1,O^2_{20},O^2_{20}\RR
\end{array}$$

\noindent
Even cocycles in grade 0 are trivial.

\noindent
Coboundary component expressions in grade 0:
$$\begin{array}{ccc}
C\LR E_{01},E_{02},E_{03},E_{04},O^1\RR&=&0\\
\vdots&\vdots&\vdots\\
C\LR O^1,O^1,O^2_{10},O^2_{10},O^2_{20}\RR&=&4t_{239}+2t_{245}
\end{array}$$
where
$$\begin{array}{ccc}
t_1&=&C\LR E_{01},E_{02},E_{05},O^1\RR\\
\vdots&\vdots&\vdots\\
t_{245}&=&C\LR E_{20},O^1,O^2_{10},O^2_{10}\RR
\end{array}$$

\noindent
Odd cocycles in grade 0:
\begin{eqnarray*}
(1)\ a^5_0&=&C\LR E_{02},E_{10},E_{11},E_{20},O^2_{01}\RR
-2C\LR E_{10},E_{11},O^1,O^2_{02},O^2_{10}\RR\\
&&
-C\LR E_{11},E_{20},O^1,O^2_{01},O^2_{01}\RR\\
&=&C\LR2yX,-Y,xX-yY,-2xY,y\RR
-2C\LR-Y,xX-yY,X Y,y^2,x\RR\\
&&-C\LR xX-yY,-2xY,X Y,y,y\RR
\end{eqnarray*}

\noindent
and also:
\begin{eqnarray*}
(1)&&
C\LR E_{01},E_{02},E_{10},E_{11},O^2_{20}\RR
+4C\LR E_{10},E_{11},O^1,O^2_{01},O^2_{11}\RR\\
&&
-2C\LR E_{10},E_{11},O^1,O^2_{02},O^2_{10}\RR\\
&&=
C\LR X,2yX,-Y,xX-yY,x^2\RR
+4C\LR-Y,xX-yY,X Y,y,x y\RR\\
&&
-2C\LR-Y,xX-yY,X Y,y^2,x\RR = a^5_0\\
\vdots&&\\
(4)&&
C\LR E_{01},E_{10},E_{11},E_{20},O^2_{02}\RR
-4C\LR E_{10},E_{11},O^1,O^2_{01},O^2_{11}\RR\\
&&
-C\LR E_{11},E_{20},O^1,O^2_{01},O^2_{01}\RR\\
&&=
C\LR X,-Y,xX-yY,-2xY,y^2\RR
-4C\LR-Y,xX-yY,X Y,y,x y\RR\\
&&
-C\LR xX-yY,-2xY,X Y,y,y\RR = a^5_0
\end{eqnarray*}

\section{Buttin vector fields and related algebras}
\noindent
The Poisson and Hamiltonian (super)algebras are very important algebras
of vector fields.
In the papers on the {\em deformation quantization}\cite{Fedosov94,Kontsevich97}
it was proven that the Poisson bracket is an {\em unique}
structure providing deformation of a commutative algebra of
differentiable functions on a manifold to a new noncommutative but
associative algebra.
In paper\cite{KornComp} we present some results about the structure of cohomology rings
for Poisson, Hamiltonian and related algebras obtained with the help of our program.

Whereas the Poisson bracket, defined on a $2n$-dimensional symplectic manifold,
has an old history, its counterpart called Buttin bracket (or odd Poisson bracket,
or antibracket) and defined on a $(n|n)$-dimensional {\em odd symplectic\/}
supermanifold\fnm{d}\fnt{d}{Such manifolds possess interesting and unusual geometrical
properties.\cite{Khudaverdian,KhudaverdianNersesian}} is a comparatively new construction.
The first example of such bracket has appeared in Schouten's paper\cite{Schouten40}
as an extension of the Lie bracket on vector fields to an bracket on
skew-symmetric contravariant (i.e., tangent) tensor fields ({\em multivectors}).
A more abstract formulation for this bracket was given by Buttin.\cite{Buttin69}
Since 1981 antibrackets are very popular in theoretical physics, because they
play the crucial role in the Batalin-Vilkovisky (BV) covariant method for quantizing
general gauge theories.\cite{BatalinVilkovisky81}
This method called the BV (or antibracket, or field-antifield)
formalism\fnm{e}\fnt{e}{The antibracket in BV formalism is an odd symplectic form on
the (infinite-dimensional) space of fields and antifields playing the role of even and
odd variables respectively, and the partial derivatives in formulas (\ref{Bb}) and
(\ref{Delta}) should be replaced by variational derivatives and the summation
by integration.}
being currently a most powerful procedure for quantizing gauge theories is applied
also in string and topological field theories.\cite{GomisParisSamuel94}

Therefore investigation of properties of Buttin and related algebras
is a problem of interest for physics. We should stress also importance
of the special subalgebras $\MR{SB}(*)$ and $\MR{SLe}(*)$, because $\Delta$-operator
(\ref{Delta}) plays an essential role in the BV formalism: so-called
{\em master equation\/} in the BV formalism is defined via this operator.

\subsection{Computations}
\noindent
Here we present some results of computations of cohomologies in the trivial
module for the algebras Buttin $\MR{B}(1)$, special Buttin $\MR{SB}(1)$
and their centerless quotients $\MR{Le}(1)$ and $\MR{SLe}(1)$.
We present also the results for odd contact algebra $\MR{M}(1)$
and its special subalgebra $\MR{SM}(1)$.
Our computations are restricted
with values for cohomology degree ($\leq 10$) and grade ($\leq 10$).
We computed also some cocycles for the case $n\geq1$. Here we mention
only the most regular of them: $a^1_{-n}=C(X_1\cdots X_n)$ for the algebras
$\MR{SB}(n)$ and $\MR{SLe}(n)$, this cocycle generates an infinite
number of its wedged powers if $n$ is even, and
$a^2_0=C(X_1,x_1)=\ldots=C(X_n,x_n)$ for the algebras
$\MR{Le}(n)$ and $\MR{SLe}(n).$ Any cocycle for the algebra $\MR{SB}(n)$
generates also new cocycles due to presence of odd center in this algebra.

We use the following grading for the variables:
$g(x_i) = 1,\ g(X_i) = -1,\ g(T) = 0.$ With this grading the algebras
$\MR{B}(n), \ \MR{Le}(n)$ and $\MR{M}(n)$ contain an internal grading
element at any $n$, i. e., there is no need to compute cohomology in
grades different from zero for this algebras. Unfortunately, there is
no good grading providing internal grading element
for the special subalgebras.\fnm{f}\fnt{f}{There are grading elements for
the algebras
$\MR{SB}(n),\ \MR{SLe}(n)$ and $\MR{SM}(n)$ for $n$ even at the grading:
$g(x_i) = 1,\  g(X_i) = -1,\ 1 \leq i \leq \F{n}{2};\ g(x_i) = -1,
\  g(X_i) = 1,\ \F{n}{2} < i \leq n;$
but in this case the space of cochains in zero grade becomes infinite-dimensional.}
One can see also that the even generating functions
correspond to the odd elements of algebra and vice versa.
In particular, the element $1$
is an {\em odd\/} central element in the algebras $\MR{B}(n)$ and
$\MR{SB}(n).$

In the below formulas all indices $i, j\geq 0$, but for the algebras
$\MR{Le}(1)$ and $\MR{SLe}(1)$ the central element $O_0=1$ should be excluded.

\subsubsection{$H^k(\MR{B}(1))$ and $H^k(\MR{Le}(1))$}
\noindent
Basis elements:
$$\begin{array}{lllll}
(1)&E_{i}&=&x^{i}X;&g(E_{i})=i+1;\\
(2)&O_{i}&=&x^{i};&g(O_{i})=i.
\end{array}$$
\noindent
Non-zero commutators:
$$\begin{array}{llll}
(1)&[E_{i}, E_{j}]&=&(i-j)E_{i+j-1}\\
(2)&[E_{i}, O_{j}]&=&-jO_{i+j-1}\\
\end{array}$$
\noindent
Generating cocycles for $H^k(\MR{B}(1)):$
\begin{eqnarray*}
a^3&=&C(X,xX,x^2X)\\
b^3&=&C(X,xX,x)-\F{1}{2}C(X,x^2X,1)
\end{eqnarray*}
The ring $H^*(\MR{B}(1))$ contains also all cocycles of the form
$a^3\W C(1)\W\cdots\W C(1)$
and  $b^3\W C(1)\W\cdots\W C(1).$

$H^k(\MR{Le}(1))$ contains only 3 non-trivial cocycles for $k \leq 20:$
the above cocycle $a^3$ and the cocycles $b^3=C(X,xX,x)$ and $a^2=C(X,x).$
The last cocycle describes the central extension of $\MR{Le}(1)$
to $\MR{B}(1).$

\subsubsection{$H^k(\MR{M}(1))$}
\noindent
Basis elements:
$$\begin{array}{lllll}
(1)&E^1_{i}&=&x^{i}X;&g(E^1_{i})=i-1;\\
(2)&E^2_{i}&=&x^{i}T;&g(E^2_{i})=i;\\
(3)&O^1_{i}&=&x^{i}T X;&g(O^1_{i})=i-1;\\
(4)&O^2_{i}&=&x^{i};&g(O^2_{i})=i.
\end{array}$$

\noindent
Non-zero commutators:
$$\begin{array}{llll}
(1)&[E^1_{i}, E^1_{j}]&=&(j-i)E^1_{i+j-1}\\
(2)&[E^1_{i}, E^2_{j}]&=&(-i+1)E^1_{i+j}+jE^2_{i+j-1}\\
(3)&[E^2_{i}, E^2_{j}]&=&(j-i)E^2_{i+j}\\
(4)&[E^1_{i}, O^1_{j}]&=&(j-i)O^1_{i+j-1}\\
(5)&[E^2_{i}, O^1_{j}]&=&(j-i+1)O^1_{i+j}\\
(6)&[E^1_{i}, O^2_{j}]&=&jO^2_{i+j-1}\\
(7)&[E^2_{i}, O^2_{j}]&=&(j-2)O^2_{i+j}\\
(8)&[O^1_{i}, O^2_{j}]&=&(-j+2)E^1_{i+j}+jE^2_{i+j-1}\\
\end{array}$$

\noindent
We have found only one non-trivial cocycle
$$\begin{array}{lll}
a^3&=&C(T,T X,x)-\F{1}{2}C(T,xT X,1)\\
&=&C(X,xX,xT)-\F{1}{2}C(T,xT X,1)\\
&=&-C(X,T,xT)-\F{1}{2}C(T,xT X,1)\\
&=&C(X,xT X,x)+\F{1}{2}C(x^2X,T X,1)+\F{1}{2}C(xT,T X,1)\\
&=&-C(xX,T X,x)-\F{1}{2}C(T,xT X,1).
\end{array}$$

\subsubsection{$H^k_g(\MR{SB}(1))$ and $H^k_g(\MR{SLe}(1))$}
\noindent
Basis elements:
$$\begin{array}{lllll}
(1)&E&=&X;&g(E)=-1;\\
(2)&O_{i}&=&x^{i};&g(O_{i})=i.
\end{array}$$

\noindent
Non-zero commutators:
$$\begin{array}{llll}
(1)&[E, O_{i}]&=&-iO_{i-1}\\
\end{array}$$

\noindent
The algebras $\MR{SB}(1)$ and $\MR{SLe}(1)$ contain the odd centers
$Z=\MR{Span}(\{1\})$ and $Z=\MR{Span}(\{x\})$ respectively.
The sets of non-trivial cocycles consist of some generating cocycles
and their consequences obtained by multiplication of these
cocycles by arbitrary wedged powers of $C(1)$ and $C(x)$ for
$\MR{SB}(1)$ and $\MR{SLe}(1)$ respectively.

Table \ref{tableS(1)} presents all generating non-trivial
cocycles in the limits for cohomology degree $k\leq 10$ and grade $g\leq 10$ for
algebras $\MR{SB}(1)$ and $\MR{SLe}(1).$
The presence of generating cocycle for $\MR{SLe}(1)$ is marked
with * in the table.
We give here the explicit expressions up to 3-cocycles for
$\MR{SB}(1)$ (The cocycles for $\MR{SLe}(1)$ can be obtained
by deleting terms with argument 1 from these expressions).
\begin{eqnarray*}
a^1_{\mbox{-}1}\!\!&=&C(X)\\
a^3_1&=&C(X,1,x^2)-C(X,x,x)\\
a^3_3&=&C(X,1,x^4)-4C(X,x,x^3)+3C(X,x^2,x^2)\\
a^3_5&=&C(X,1,x^6)-6C(X,x,x^5)
+15C(X,x^2,x^4)-10C(X,x^3,x^3)\\
a^3_7&=&C(X,1,x^8)-8C(X,x,x^7)+28C(X,x^2,x^6)
-56C(X,x^3,x^5)\\&&
+35C(X,x^4,x^4)\\
a^3_9
&=&C(X,1,x^{10})-10C(X,x,x^9)
+45C(X,x^2,x^8)-120C(X,x^3,x^7)\\&&
+210C(X,x^4,x^6)-126C(X,x^5,x^5)
\end{eqnarray*}
\begin{table}[h!]
\caption{Generating cocycles for $H^k_g(\MR{SB}(1))$
 and $H^k_g(\MR{SLe}(1))$}
\begin{center}\label{tableS(1)}
\begin{tabular}{|c|l|l|l|l|l|l|l|l|l|l|l|l|}
\hline
$k\backslash g$&-1&0&1&2&3&4&5&6&7&8&9&10
\Overline
1&$a^{1~~*}_{-1}$&&&&&&&&&&&
\Interline
2&&&&&&&&&&&&
\Interline
3&&&$a^{3~*}_1$&&$a^{3~*}_3$&&$a^{3~*}_5$&&$a^{3~*}_7$&&$a^{3~*}_9$&
\Interline
4&&&&$a^4_2$&&$a^4_4$&$a^{4~*}_5$&$a^4_6$&$a^{4~*}_7$&$a^{4~*}_8$&$a^{4~*}_9$&$a^{4~~*}_{10}$
\Innerline
&&&&&&&&&&$b^4_8$&&$b^4_{10}$
\Interline
5&&&&&$a^5_3$&&$a^5_5$&$a^5_6$&$a^{5~*}_7$&$a^5_8$&$a^{5~*}_9$&$a^{5~~*}_{10}$
\Innerline
 &&&&&       &&       &       &$b^5_7$    &       &$b^5_9$    &$b^5_{10}$
\Innerline
 &&&&&       &&       &       &           &       &$c^5_9$    &
\Interline
6&&&&&&$a^6_4$&&$a^6_6$&$a^6_7$&$a^6_8$&$a^{6~*}_9$&$a^6_{10}$
\Innerline
 &&&&&&       &&       &       &$b^6_8$&$b^6_9$    &$b^6_{10}$
\Innerline
 &&&&&&&&&&&&$c^6_{10}$
\Interline
7&&&&&&&$a^7_5$&&$a^7_7$&$a^7_8$&$a^7_9$&$a^7_{10}$
\Innerline
 &&&&&&&       &&       &       &$b^7_9$&$b^7_{10}$
\Interline
8&&&&&&&&$a^8_6$&&$a^8_8$&$a^8_9$&$a^8_{10}$
\Innerline
 &&&&&&&&&&&&$b^8_{10}$
\Interline
9&&&&&&&&&$a^9_7$&&$a^9_9$&$a^9_{10}$
\Interline
10&&&&&&&&&&$a^{10}_8$&&$a^{10}_{10}$ \\[-3mm]
  &&&&&&&&&&&& \\
\hline
\end{tabular}
\end{center}
\end{table}

\subsubsection{$H^k_g(\MR{SM}(1))$}

\noindent
Basis elements:
$$\begin{array}{lllll}
(1)&E_{i}&=&ix^{i-1}T+(-i+2)x^{i}X;&g(E_{i})=i-1;\\
(2)&O^1&=&T X;&g(O^1)=-1;\\
(3)&O^2_{i}&=&x^{i};&g(O^2_{i})=i.
\end{array}$$

\noindent
Non-zero commutators:
$$\begin{array}{llll}
(1)&[E_{i}, E_{j}]&=&(2j-2i)E_{i+j-1}\\
(2)&[E_{i}, O^2_{j}]&=&(2j-2i)O^2_{i+j-1}\\
(3)&[O^1, O^2_{i}]&=&E_{i}\\
\end{array}$$
Our computations in the limits $k,g\leq10$ show that
the only non-trivial cocycles take the form $a^k_{k\mbox{-}2}$
excepting the case $k=2$. We give below the explicit expressions
for the first three these cocycles. The cocycle $a^1_{\mbox{-}1}$
generates an infinite set of cocycles
$a^k_{\mbox{-}k} = a^1_{\mbox{-}1}\W\cdots\W a^1_{\mbox{-}1}
= C(TX,\ldots,TX).$

\begin{eqnarray*}
a^1_{\mbox{-}1}\!\!
 &=&C(T X)\\
a^3_1
 &=&C(T+xX,2xT,1)-2C(T X,x,x)\\
 &=&C(2X,T+xX,x^2)-2C(T X,x,x)\\
 &=&-C(2X,2xT,x)+2C(T X,1,x^2)+2C(T X,x,x)\\
a^4_2
 &=&
  C(T+xX,2xT,1,x)-\F{1}{6}C(T+xX,3x^2T-x^3X,1,1)-\F{4}{3}C(T X,x,x,x)\\
 &=&
  C(2X,T+xX,x,x^2)-\F{1}{6}C(T+xX,3x^2T-x^3X,1,1)-\F{4}{3}C(T X,x,x,x)\\
 &=&
  -\F{1}{2}C(2X,2xT,x,x)-\F{1}{6}C(T+xX,3x^2T-x^3X,1,1)+2C(TX,1,x,x^2)\\
  & &+\F{2}{3}C(TX,x,x,x)
\end{eqnarray*}

\section{Conclusion}
\noindent
Mathematicians regard the problem of
computation of cohomology as solved if they construct the full
cohomology ring. The structure of such rings may be rather
complicated including many non-trivial relations between the cocycles
in contrast to the examples in this paper where the rings are
commutative. To get a clear idea about the structure of cohomology
ring one should compute usually the cocycles up to
degrees high enough. Unfortunately the computation of cohomology
is a typical
problem with the combinatorial explosion. Nevertheless, some results
can be obtained with the help of computer having an efficient
enough program. On the other hand, physicists are interested mainly
in the
second cohomologies describing the central extensions and deformations.
Such cohomologies can be computed rather easily even for large algebras.
Some essential possibilities remain for increasing the efficiency
of the program and we hope to implement the corresponding
improvements in future.

\nonumsection{Acknowledgment}
\noindent
I would like to thank D. Leites for initiating this work
and helpful communications.
I am also grateful to V. Gerdt and O. Khudaverdian for fruitful
discussions and useful advises.
This work was supported in part by INTAS project No. 96-0842 and
RFBR project No. 98-01-00101.

\newpage
\nonumsection{References}


\vspace*{-0.25cm}
\begin{thebibliography}{000}
%1.
\bibitem{Fuks}
D.B. Fuks, {\bibit Cohomology of Infinite Dimensional Lie Algebras\/}
          (Consultants Bureau, New York, 1987).
%2.
\bibitem{GL}
P. Grozman and D. Leites, in {\bibit The second International
Mathematica symposium}, ed. V. Ker\"anen (Rovaniemi, 1997), p. 185.
%3.
\bibitem{Fil}
V.T. Filippov, {\bibit Sibirskii Math. J.}, {\bibbf 24}, 126 (1985).
(in Russian)
%4.
\bibitem{Takh}
L. Takhtajan, {\bibit Commun. Math. Phys.}, {\bibbf 160}, 295 (1994).
%5.
\bibitem{D'Hoker}
E. D'Hoker, {\bibit Nucl. Phys.}, {\bibbf B451}, 725 (1998).
%6.
\bibitem{Leites}
D. Leites, in {\bibit Modern Problems of Mathematics.
 Recent developments}, {\bibbf 25}, VINITI, Moscow, 1984, p. 3
(in Russian; English translation in JOSMAR {\bibbf 30(6)}, 1985, p. 2481)
%7
\bibitem{Pommaret}
J.-F. Pommaret, {\bibit Partial Differential Equations and Group Theory\/}
          (Kluwer, Dordrecht, 1994).
%8.
\bibitem{LP}
D. Leites and G. Post, in {\bibit Computers and Mathematics},
eds. E. Kaltofen and S.M. Watt (Springer, NY ea, 1989), p. 73.
%9.
\bibitem{PH}
G. Post and N. von Hijligenberg, ``Calculation of Lie algebra cohomology
by computer", Memo\# 833, Faculty of Appl. Math. Univ. Twente, 1989;
id. ibid. \#928, 1991.
%10.
\bibitem{FL}
D. Fuchs and D. Leites, {\bibit C.r. Acad. Bulg. Sci.}, {\bibbf 37},
1595 (1984).
%11.
\bibitem{KornProg}
V. V. Kornyak, {\bibit Zapiski nauchnyh seminarov POMI}, {\bibbf 258}, 148
(St.Petersburg, 1999).
%12.
\bibitem{Fedosov94}
B. Fedosov, {\bibit J. Diff. Geom.}, {\bibbf 40}, 213 (1994).
%13.
\bibitem{Kontsevich97}
M. Kontsevich, ``Deformation Quantization of Poisson Manifolds I",
q-alg/9709040.
%14.
\bibitem{KornComp}
V. V. Kornyak, in {\bibit Computer Algebra in Scientific Computing}, eds. V.G.Ganzha,
 E.W.Mayr and E.V.Vorozhtsov, (Springer, 1999), p. 241;\\
``Cohomology of Lie Superalgebras of Hamiltonian Vector Fields:
Computer Analysis", math.SC/9906046.
%15.
\bibitem{Khudaverdian}
O.M. Khudaverdian, {\bibit J. Math. Phys.}, {\bibbf 32}, 1934 (1991);
 {\bibit Commun. Math. Phys.}, {\bibbf 198}, 591 (1998); math.DG/9909117.
%16.
\bibitem{KhudaverdianNersesian}
O.M. Khudaverdian and A.P. Nersesian, {\bibit Mod. Phys. Lett.}, {\bibbf A8}, 2377 (1993);
 {\bibit J. Math. Phys.}, {\bibbf 37}, 3713 (1996).
%17.
\bibitem{Schouten40}
J.A. Schouten, {\bibit Proc. Kon. Ned. Akad. Wet. Amsterdam.} {\bibbf 43}, 449 (1940).
%18.
\bibitem{Buttin69}
C. Buttin, {\bibit C. R. Acad. Sci. Paris, Ser. A-B}, {\bibbf 269}, A87 (1969).
%19.
\bibitem{BatalinVilkovisky81}
I.A. Batalin and G.A. Vilkovisky, {\bibit Phys. Lett.}, {\bibbf 102B}, 27 (1981).
%20.
\bibitem{GomisParisSamuel94}
J. Gomis, J. Paris, S. Samuel, {\bibit Phys. Rept.}, {\bibbf 259},
1 (1995); hep-th/9412228.
\end{thebibliography}
\end{document}